\documentclass[letterpaper, 10 pt, conference]{ieeeconf}  % Comment this line out
\IEEEoverridecommandlockouts                              % This command is only
\usepackage[para,online,flushleft]{threeparttable}
\usepackage{algorithm}
\usepackage[noend]{algorithmic}
\usepackage{booktabs}
\usepackage{multirow}
\usepackage{adjustbox}
\usepackage{graphicx}
\usepackage{siunitx}
\usepackage[utf8]{inputenc}
\usepackage[T1]{fontenc}
\usepackage{amsfonts}
\usepackage{gensymb}

\usepackage{amsthm}

\usepackage{epsfig} % for postscript graphics files
\usepackage{amsmath} % assumes amsmath package installed
\usepackage{amssymb}  % assumes amsmath package installed
\usepackage{siunitx}
\usepackage{subcaption}
\usepackage{caption}
\usepackage{color}
\usepackage[table]{xcolor}
\usepackage{soul}
\usepackage{theoremref}
\usepackage[noadjust]{cite}
\usepackage{amsthm}
\usepackage{mathtools}
\DeclarePairedDelimiter{\ceil}{\lceil}{\rceil}

\usepackage{hyperref}
\hypersetup{
    colorlinks=true,
    linkcolor=blue,
    filecolor=magenta,      
    urlcolor=blue,
    pdftitle={Overleaf Example},
    pdfpagemode=FullScreen,
    }

\newtheorem*{theorem*}{Theorem}

\title{\LARGE\bf Benchmark Results for Bookshelf Organization Problem as Mixed Integer Nonlinear Program with Mode Switch and Collision Avoidance}
\author{Xuan Lin$^{1}$, Gabriel I.~Fernandez$^{1}$, and Dennis W.~Hong$^{1}$% <-this % stops a space
\thanks{$^{1}$X. Lin, Gabriel I.~Fernandez, and Dennis W.~Hong are with the Robotics and Mechanisms Laboratory, Department of Mechanical and Aerospace Engineering, University of California, Los Angeles, CA 90095, USA.
        {\tt\small \{maynight,gabriel808,dennishong\}@ucla.edu}}
}

\begin{document}
\maketitle
\thispagestyle{empty}
\pagestyle{empty}
%%%%%%%%%%%%%%%%%%%%%%%%%%%%%%%%%%%%%%%%%%%%%%%%%%%%%%%%%%%%%%%%%%%%%%%%%%%%%%%%
%%note: must be manage risk in a "quantitative" way! Model and environmental uncertainty are coupled.
%% reduce to simpler deterministic optimization when I assume the noise follows gauusian

\begin{abstract}
Mixed integer convex and nonlinear programs, MICP and MINLP, are expressive but require long solving times. Recent work that combines data-driven methods on solver heuristics has shown potential to overcome this issue allowing for applications on larger scale practical problems. To solve mixed-integer bilinear programs online with data-driven methods, several formulations exist including mathematical programming with complementary constraints (MPCC), mixed-integer programming (MIP). In this work, we benchmark the performances of those data-driven schemes on a bookshelf organization problem that has discrete mode switch and collision avoidance constraints. The success rate, optimal cost and solving time are compared along with non-data-driven methods. Our proposed methods are demonstrated as a high level planner for a robotic arm for the bookshelf problem.

\end{abstract}
% 
% 
% \begin{IEEEkeywords} Legged Robots, Motion and Path Planning,  Optimization and Optimal Control
%  \end{IEEEkeywords}

%%%%%%%%%%%%%%%%%%%%%%%%%%%%%%%%%%%%%%%%%%%%%%%%%%%%%%%%%%%%%%%%%%%%%%%%%%%%%%%%
\section{Introduction}
\label{Sec:introduction}
Optimization-based methods are useful tools for solving robotic motion planning problems. Typical approaches such as mixed-integer convex programs (MICPs) \cite{deits2014footstep, lin2019optimization}, nonlinear or nonconvex programs (NLPs) \cite{dai2014whole,winkler2018gait,shirai2020risk} and mixed-integer NLPs (MINLPs) \cite{soler2011route} offer powerful tools to formulate these problems. However, each has its own drawbacks. NLPs tend to converge to local optimal solutions. In practice, local optimal solutions can have bad properties, such as inconsistent behavior as they depend on initial guesses. Mixed-integer programs (MIPs) explicitly treat the discrete variables. Branch-and-bound is usually used to solve MIPs \cite{boyd2007branch}. MIP solvers seek global optimal solutions with more consistent behavior than NLP solvers. For small-scale problems, these algorithms usually find optimal solutions within a reasonable time \cite{lin2019optimization,tordesillas2019faster}. On the contrary, MIPs can require impractically long solving times for problems with a large number of integer variables \cite{lin2021designing}. MINLPs incorporate both integer variables and nonlinear constraints, hence, very expressive. Unfortunately, we lack efficient algorithms to tackle MINLPs. As a result, it is difficult to implement most of optimization schemes online.

\begin{figure}[!t]
		\centering
		\hspace*{-0.45cm}
		\includegraphics[scale=0.28]{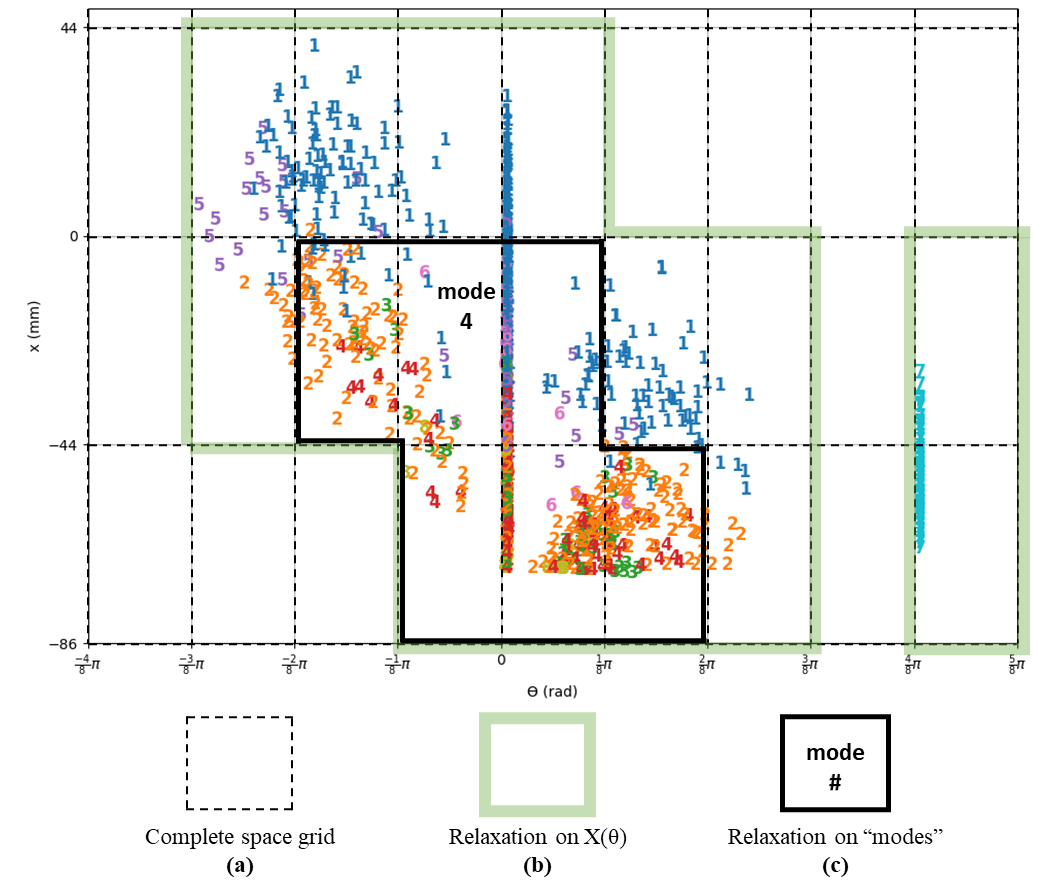}
		\caption {Angle and x position of solutions to the bookshelf problem's furthest left book are depicted by colored numbers indicating cluster membership. Each box on the grid corresponds to a different relaxation scheme using integer variables: (a) full relaxation of the whole space, (b) relaxation over the space with data, (c) relaxation for specific modes from clusters. For visualization purposes only mode 4 region is explicitly displayed. Notice how different modes can reduce the amount of integer variables needed.}
		\label{Fig:Figure_1}
\end{figure}

Recently, researchers have started to investigate machine learning methods to gather problem specific heuristics and speed up the optimization solving process. Standard algorithms to solve MIPs such as branch-and-bound and cutting plane methods rely on heuristics to quickly remove infeasible regions. Learning methods can be used to acquire better heuristics. For example, \cite{nair2020solving} used graph neural networks to learn heuristics. \cite{tang2020reinforcement} used reinforcement learning to discover efficient cutting planes. On the other hand, data can be collected to learn and solve specific problems \cite{zhu2019fast,cauligi2021coco}. In \cite{cauligi2021coco}, the authors proposed CoCo which collects problem inputs and solved integer variable strategies offline and trains a neural network. Effectively, this becomes a classification problem, where each strategy has a unique label. For online solving, the neural network will propose candidate solutions, reducing it to a convex program. The results show that CoCo can solve MIPs with around 50 integer variables within a second. However, CoCo has relatively simple neural-network structure. Since the number of possible strategies equals to $2^\lambda$, where $\lambda$ is the number of integer variables, the number of unique strategies tend to be close to the total amount of data. This induces overfitting.

Collision avoidance constraints typically appear in multi-agent motion planning problems. Single-agent collision avoidance with the stationary environment is relatively straightforward as the collision-free regions can be precomputed such that the moving agent stays within them. Polyhedrons are typical models for collision-free environment \cite{deits2015computing}. If moving agents have simple shapes such as circles, collision-free constraints can simply enforce that distance between two agents are larger than a known value \cite{dai2017distributed}. However, if moving agents have complex shapes, collision avoidance becomes more challenging. Naive methods include enforcing collision avoidance between a large number of points sampled within each agent. This approach is obviously computationally expensive. Methods based on KKT conditions are proposed \cite{raghunathan2022pyrobocop} to enforce collision avoidance between convex shapes. This method does not allow the agent to make contact without sacrificing precision. In this paper, we propose a novel approach using separating planes to enforce collision avoidance between convex agents but still allow them to make precise contacts. If the agent has a non-convex shape, it can be decomposed into convex shapes such that this approach can still be used.

For many multi-agent optimization problems, agents can operate under various modes encoded with discrete variables. The discrete mode switching with collision avoidance constraints leads to a mixed-integer bilinear formulation which is challenging to solve. Leveraging on pre-solved data is one approach to solve these formulations fast and reliably. To get a comprehensive understanding of the performance of various data-driven methods on such problems, we introduce the bookshelf organization problem for a benchmark. Given a bookshelf with several books on top, an additional book needs to be placed on the shelf with minimal disturbance to the existing books. The bookshelf problem works well as a good benchmark because it: 1) is an MINLP that can be converted to an MICP problem with hundreds of integer variables, 2) can easily be scaled to push algorithms to their limit, and 3) has practical significance where data can reasonably be gathered, such as in the logistics industry.
 
To summarize, our contributions are as follows:
\begin{enumerate}
    \item Formulate the bookshelf organization problem as an MINLP and solve it within seconds. 
    \item Benchmark the performance of different data-driven formulations on the bookshelf organization problem.
    \item Implement the solutions on hardware to stow items into the shelf.
\end{enumerate}

\section{Related Works}
\label{Sec:related_work}
\subsection{Motion Library}
Motion library is a typical method to leverage offline computing power. A large number of controllers are pre-computed offline and stored inside a library. Online, the controllers are selected based on the robot state in the nearest neighbor manner. This method can provide more robust and global policies faster than online solving techniques such as dynamic programming \cite{stolle2006policies}. Motion libraries have been used to synthesize controllers for bipedal walking \cite{liu2013biped}, and simulated multi-link swimming \cite{tassa2007receding}. However, this method does not adapt well to situations not previously seen in the library.

\subsection{Parametric Programming}
Parametric programming is a technique to build a function that maps to the optimization solutions from varying parameters of the optimization problem \cite{fiacco1983introduction, pistikopoulos2002line}. Previous works has investigated parametric programming on linear programming \cite{gal2010postoptimal}, quadratic programming \cite{bemporad2002explicit}, mixed-integer nonlinear programming \cite{dua1999algorithms}. As an implementation of parametric programming on controller design, explicit MPC \cite{bemporad2002explicit, tondel2003algorithm} tries to solve the model-predictive control problem offline, and stores active sets. When computed online, the problem parameter can be used to retrieve the active sets. \cite{bemporad2002explicit} solved a constrained linear quadratic regulator problem with explicit MPC and proved that the active sets are polyhedrons. Therefore, the parameter space can be partitioned using an algorithm proposed by \cite{dua2000algorithm}. However, when the problem is non-convex, the active sets do not in general form polyhedrons. Therefore, it is significantly more complex to compute critical regions offline.

\subsection{Learning Problem-solution Mapping}
The basic ideas behind parametric optimization is to construct a mapping from the descriptor of the problem to the solutions of the problem. For simple problems, this may be done analytically with parametric optimization techiques. For more challenging problems, it makes more sense to use learning approaches. For relatively smaller datasets, previous works \cite{hauser2016learning, tang2019data, zhu2019fast} directly store the data points and pick out warm-start using non-parametric learning such as K-nearest neighbor (KNN), and solve the online formulation. Effectively, this approach adds an online adaptation step to the motion library method, which turns out to work well. On the other hand, modern learning techniques such as neural-networks can learn an embedding of a larger set of parameters that maps to the solutions \cite{cauligi2021coco}. One advantage of using neural-network methods is the capability to deal with out-of-distribution situations that are not included in the training set \cite{power2022variational}. 

\subsection{Mixed-integer Programming}
To allow online adaptation, one needs to solve an online optimization formulation in real-time. This online formulation can have several forms. Typical robotics problem includes discrete contact and discrete operating modes. There are generally two ways to include discrete variables in the optimization formulation. Mixed-integer programs include continuous variables and explicit discrete variables. The definition of mixed-integer variables is that if the discrete variables are relaxed into continuous ones, the problem becomes convex \cite{dai2019global}. Typical mixed-integer programming solvers use branch-and-bound methods \cite{boyd2007branch}, cutting plane methods \cite{marchand2003cutting}, and rely on heuristics to choose the branch to explore. Despite the worst case of solving time, a large portion of problems only requires exploring a small portion of the search tree \cite{williams2013model}. Mixed-integer programs have been implemented for online motion planning such as \cite{tordesillas2019faster}.

\subsection{Mathematical Programming with Complementary Constraints}
On the other hand, mathematical programs with complementary constraints (MPCC) models discrete modes through continuous variables with complementary constraints. Complementary constraints enforce a pair of variables such that if one of them is non-zero, the other one should be zero. This constraint is traditionally hard to solve and is sensitive to initial guesses. Algorithms such as time-stepping \cite{anitescu1997formulating}, pivoting \cite{drumwright2015rapidly}, central path methods \cite{kojima1991unified} are proposed to resolve complementarity. In the robotics community, complementary constraints are typically used to optimize over gaits for trajectory optimization \cite{posa2014direct, zhang2021transition} or control with implicit contacts \cite{cleac2021fast} where contact forces and distance to the ground are complementary with each other. On the other hand, complementary constraints can also be used to model binary variables \cite{park2018semidefinite}.

\section{Bookshelf Organization Problem Setup}
\label{Sec:problem_setup}
\begin{figure*}[!t]
		\centering
		\hspace*{-0.25cm}
		\includegraphics[scale=0.52]{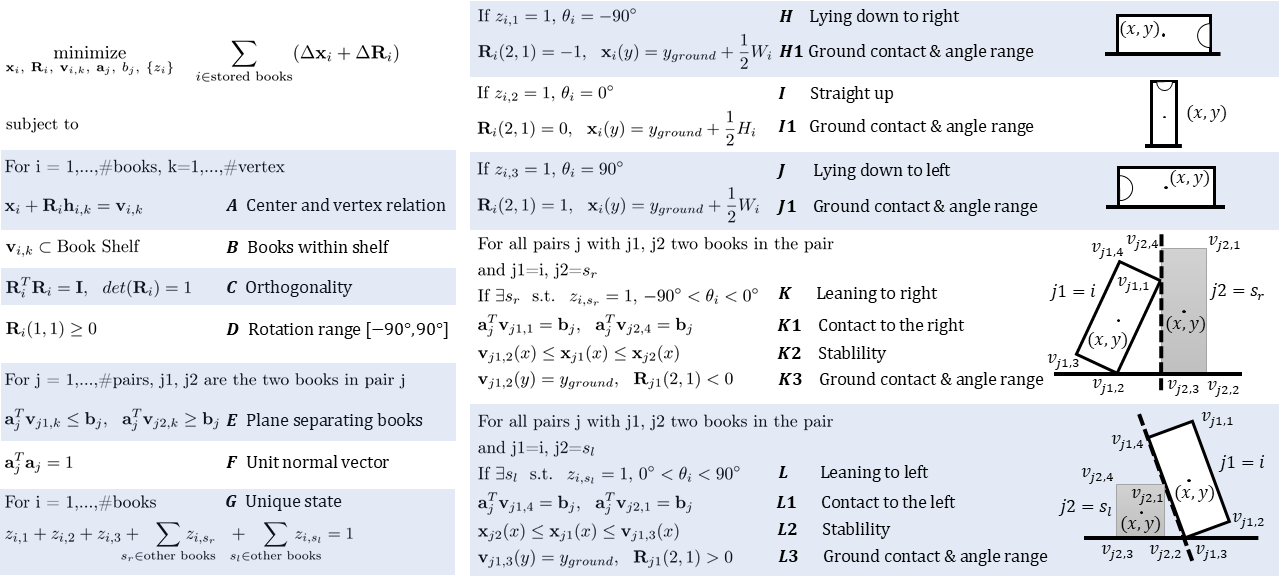}
		\caption {Complete formulation of the bookshelf organization problem.}
		\label{Fig:Complete_formulation}
\end{figure*}

Assume a 2D bookshelf with limited width $W$ and height $H$ contains rectangular books where book $i$ has width $W_{i}$ and height $H_{i}$ for $i=1,...,N-1$. A new book, $i=N$, is to be inserted into the shelf. The bookshelf contains enough books in various orientations, where in order to insert book $N$, the other $N-1$ books may need to be moved, i.e., optimize for minimal movement of $N-1$ books. This paper focuses on inserting one book using a single robot motion during which the book is always grasped. The problem settings can be extended to a sequence of motions with re-grasping. This problem is useful in the logistics industry.
%, such as robots filling shelves.

% In this section, we introduce the formulation of the bookshelf organization problem that we focus on in this paper.

% Assume a squared 2D bookshelf with width $W$ and height $H$ contains a few books of squared shape. Book $i$ has width $W_{i}$ and height $H_{i}$, $i=1,...,N$. A new book $i=N$ is to be inserted into the shelf. In many cases, the bookshelf has limited space and stored books need to be moved to create room for the new book. We want to optimize for minimal movement of the stored books for the new book. This problem has practical significance such as in logistic industry. For example, with the solved bookshelf plan, we can plan a trajectory and implement robot manipulators to realize the insert motion.

Fig. \ref{Fig:Complete_formulation} shows the constraints, variables and objective function for the bookshelf problem. The variables that characterize book $i$ are: position $\textbf{x}_{i} = [x_{i}, y_{i}]$ and angle $\theta_{i}$ about its centroid. $\theta_{i}=0$ when a book stands upright. The rotation matrix is: $\textbf{R}_{i} = [cos(\theta_{i}), \ -sin(\theta_{i}); \ sin(\theta_{i}), \ cos(\theta_{i})]$. Let the 4 vertices of book $i$ be $\textbf{v}_{i,k}$, $k=1,2,3,4$. The constraint \textbf{A} in Fig. \ref{Fig:Complete_formulation} shows the linear relationship between $\textbf{x}_{i}$ and $\textbf{v}_{i,k}$, where $\textbf{h}_{i,k}$ is the constant offset vector from its centroid to vertices. Constraint \textbf{B} enforces that all vertices of all books stay within the bookshelf, a linear constraint. Constraint \textbf{C} enforces the orthogonality of the rotation matrix, a bilinear (non-convex) constraint. Constraint \textbf{D} enforces that the angle $\theta_{i}$ stays within $[-90\degree, 90\degree]$, storing books right side up.

To ensure that the final book positions and orientations do not overlap with each other, separating plane constraints are enforced. For convex shapes, the two shapes do not overlap with each other if and only if there exists a separating hyperplane $\textbf{a}^{T}\textbf{x}=b$ in between \cite{boyd2004convex}. That is, for any point $\textbf{p}_{1}$ inside shape 1 then $\textbf{a}^{T}\textbf{p}_{1} \leq b$, and for any point $\textbf{p}_{2}$ inside shape 2 then $\textbf{a}^{T}\textbf{p}_{2} \geq b$. This is represented by constraint \textbf{E}. Constraint \textbf{F} enforces $\textbf{a}$ to be a normal vector. Both \textbf{E} and \textbf{F} are bilinear constraints.

% To ensure that the final book positions and orientations do not overlap with each other, we enforce separating plane constraints. For convex shapes, the two shapes do not overlap with each other if and only if there exists a separating hyperplane $\textbf{a}^{T}\textbf{x}=b$ in between \what{reference needed}. That is, for any point $\textbf{p}_{1}$ inside shape 1 then $\textbf{a}^{T}\textbf{p}_{1} \leq b$, and for any point $\textbf{p}_{2}$ inside shape 2 then $\textbf{a}^{T}\textbf{p}_{1} \geq b$. This is represented by constraint E. Constraint F enforces $\textbf{a}$ to be a normal vector. Both E and F are bilinear constraints.

% To ensure that the solved book shapes does not overlap with each other, we enforce separating plane constraints. For convex shapes, the two shapes does not overlap with each other if and only if there exist a separating hyperplane $\textbf{a}^{T}\textbf{x}=b$ in between \what{reference needed}. That is, for any point $\textbf{p}_{1}$ inside shape 1, we have $\textbf{a}^{T}\textbf{p}_{1} \leq b$, and for any point $\textbf{p}_{2}$ inside shape 2, we have $\textbf{a}^{T}\textbf{p}_{1} \geq b$. This is represented by constraint E. Constraint F enforces $\textbf{a}$ to be a normal vector. Note both E and F is another bilinear constraint.

Finally, we need to assign a state to each book $i$. For each book, it can be standing straight up, laying down on its left or right, or leaning towards left or right against some other book, as shown in the far right column in Fig. \ref{Fig:Complete_formulation}. For each book $i$, we assign a set of integer variables ${z_{i}}$. If book $i$ stands upright ($z_{i,2}=1$) or lays flat on its left ($z_{i,1}=1$) or right ($z_{i,3}=1$), constraints \textbf{I1} or \textbf{J1} or \textbf{H1} are enforced, respectively. If book $i$ leans against another book (or left/right wall which can be treated as static books) on the left or right, constraints in \textbf{K} and \textbf{L} are enforced, respectively. To this end checks need to indicate the contact between books. By looking at the right column in Fig. \ref{Fig:Complete_formulation}, we can reasonably assume that the separating plane $\textbf{a}^{T}\textbf{x}=b$ always crosses vertex 1 of the book on the left and vertex 4 of the book on the right. This is represented by bilinear constraints, \textbf{K1} and \textbf{L1}. In addition, the books need to remain stable given gravity. Constraints \textbf{K2} and \textbf{L2} enforce that a book is stable if its $x$ position stays between the supporting point of itself (vertex 2 if leaning rightward and vertex 3 if leaning leftward) and the $x$ position of the book that it is leaning onto. Lastly, constraint \textbf{K3} and \textbf{L3} enforce that the books have contact with the \emph{ground}. For practical reasons, we assume that books cannot stack onto each other, i.e, each book has to touch the \emph{ground} of bookshelf at at least one point. We note that constraints in \textbf{H}, \textbf{I}, \textbf{J}, \textbf{K}, and \textbf{L} can be easily formulated as MICP constraints using big-M formulation \cite{vielma2015mixed}, such that they are enforced only if the associated integer variable $z=1$. Also, it can easily be extended to allow stacking for our problem. Any contact conditions between pairs of books may also be added into this problem as long as it can be formulated as mixed-integer convex constraint. Overall, this is a problem with integer variables, ${z_{i}}$, and non-convex constraints \textbf{C}, \textbf{E}, \textbf{F}, \textbf{K1}, and \textbf{L1}, hence, an MINLP problem.

Practically, this problem presents challenges for retrieving high quality solutions. If robots were used to store books, the permissible solving time is several seconds, and less optimal solutions means longer realization times. The video of our hardware implementation demonstrates non-optimal insertion cases that requires multiple additional robot motions that dramatically increase the chance of failure. There are several potential approaches to resolving this issue: fix one set of nonlinear variables and solve MICP \cite{posa2015stability}, convert the nonlinear constraints into piece-wise linear constraints and formulate them into an MICP \cite{deits2014footstep}, or directly apply MINLP solvers such as BONMIN \cite{bonami2007bonmin}. As expected, these approaches struggle to meet the requirements. 

\section{Learning Algorithm} 
\label{Sec:learning_algorithm}
Assume that we are given a set of problems parametrized by $\Theta$ that is drawn from a distribution $D(\Theta)$. For each $\Theta$, we seek a solution to the optimization problem:

\begin{equation}
\begin{aligned}
& \underset{\textbf{x}, \ \textbf{z}}{\text{minimize}} \ f_{obj}(\textbf{x}, \textbf{z}; \Theta) \\
\text{s. t.} \ \ & f_{i}(\textbf{x}, \textbf{z}; \Theta) \leq 0, \ \ i = 1,...,m_{f} \\
& b_{j}(\textbf{x}, \textbf{z}; \Theta) \leq 0, \ \ j = 1,...,m_{b}
\label{Eqn:General_formulation}
\end{aligned}
\end{equation}

Where $\textbf{x}$ denotes continuous variables and $\textbf{z}$ binary variables with $z_{i} \in \{0, 1\}$ for $i=1,...,dim(\textbf{z})$. Constraints $f_{i}$ are mixed-integer convex, meaning if the binary variables $\textbf{z}$ are relaxed into continuous variables $\textbf{z} \in [0, 1]$, $f_{i}$ becomes convex. Constraints $b_{j}$ are mixed-integer bilinear, meaning that relaxing the binary variables gives bilinear constraints. Without loss of generality, $\textbf{x}$ and $\textbf{z}$ are assumed to be involved in each constraint. We omit equality constraints in \eqref{Eqn:General_formulation} as they can be turned into two inequality constraints from opposite directions. 

In general, there are two directions to convert a mixed-integer bilinear program to either a mixed-integer linear program or a nonlinear program. To convert to a mixed-integer linear program, one needs to convert non-convex constraints to mixed-integer linear constraints. For example, we can convert trigonometry constraints into piece-wise linear constraints \cite{deits2014footstep}, or convert bilinear constraints into mixed-integer envelope constraints \cite{dai2019global}. However, mixed-integer envelope constraint tend to significantly increase the solving time. In \cite{lin2021designing}, the trajectory planner is based on mixed-integer formulation with envelopes to approximate the bilinear moment constraints. To walk a few steps upstairs, the solver can take several hours to find a feasible solution. The other approach is to convert the binary variables into continuous variables with complementary constraints. This approach requires the solver to deal with a potentially large number of complementary constraints known to be computationally difficult. Directly solving the NLP with complementary constraints without an informed initial guess result in high chance of infeasibility, unless special treatment such as \cite{anitescu1997formulating} are used. Those methods, however, result in slower solving speed from our benchmark results.

For both directions mentioned above, incorporating pre-solved data to learn a good warm-start online can improve the solving speed or rate of feasibility, as we will show in our experiments. If the learning agent with the problem description $\Theta$ can partially provide feasible integer variables, the size of the mixed-integer optimization problem can be reduced and hence solved faster. If the learning agent can provide a complete list of integer variables, the problem reduces to a convex optimization permitting convex solvers such as OSQP \cite{stellato2020osqp}. On the other hand, if the learning agent provides a good initial guess for the variables involved in complementary constraints, the feasibility to solve the problem significantly increases.

To better understand the performances of the approaches mentioned above, we compared their performances on the bookshelf problem. In addition, we have also tested the performance of non-data-driven methods, including manually designed heuristics, and nonlinear ADMM \cite{risk_aware} as benchmarks. We first collect the initial dataset. This is done by sampling the problem parameters $\Theta$, then solving the samples with a non-data-driven approach. Either manually designed heuristics or nonlinear ADMM can be used. We can then use the initial dataset as warm-start points to solve the new problem instances and enlarge the dataset. With warm-start points, data-driven methods can be used, and the new dataset has better optimal costs. This procedure iterates until we arrive at the final dataset. The final dataset is used to learn good problem-solution mapping for warm-start online. The online formulation can be MPCC, MIP, or a convex formulation. In the following sections, we describe the implementation details for each online formulation. 

\begin{table*}[]
\caption{Benchmark Results}
\resizebox{\textwidth}{!}{\begin{tabular}{@{}l|c|c|c|c|cc|cc|cc|c|c|c|c@{}}
\toprule
                      & ADMM - Knitro        & ADMM - IPOPT         & MIP - 50                                                                              & MIP - 20                                                                              & \begin{tabular}[c]{@{}c@{}}MPCC\_default\\ - Knitro\end{tabular} & \begin{tabular}[c]{@{}c@{}}MPCC\_default\\ - IPOPT\end{tabular} & \begin{tabular}[c]{@{}c@{}}MPCC\_Manual\\ - Knitro\end{tabular} & \begin{tabular}[c]{@{}c@{}}MPCC\_Manual\\ - IPOPT\end{tabular} & \begin{tabular}[c]{@{}c@{}}MPCC\_KNN3\\  - Knitro (15000data)\end{tabular} & \begin{tabular}[c]{@{}c@{}}MPCC\_KNN3\\  - IPOPT (15000 data)\end{tabular} & \begin{tabular}[c]{@{}c@{}}MPCC\_KNN3 (Rev)\\ - Knitro (36000 data)\end{tabular}                          & \begin{tabular}[c]{@{}c@{}}Convex - 200\\ (15000 data)\end{tabular}                     & \begin{tabular}[c]{@{}c@{}}Convex - 600\\ (15000 data)\end{tabular}                     & \begin{tabular}[c]{@{}c@{}}Convex - 200\\ (78000 data)\end{tabular}                 \\ \midrule
Note                  & \multicolumn{1}{l}{} & \multicolumn{1}{l|}{} & \multicolumn{1}{l|}{\begin{tabular}[c]{@{}l@{}}Each cluster\\ 50 samples\end{tabular}} & \multicolumn{1}{l|}{\begin{tabular}[c]{@{}l@{}}Each cluster\\ 20 samples\end{tabular}} & \multicolumn{2}{l|}{\begin{tabular}[c]{@{}l@{}}Using default initial guess from \\ CasADi\end{tabular}}                             & \multicolumn{2}{l|}{Using manually selected 1 initial guess}                                                                      & \multicolumn{2}{l|}{KNN selects top 3 guesses from 15000 data}                                                                                           & \multicolumn{1}{l|}{\begin{tabular}[c]{@{}l@{}}Same to left but select\\ the worst 3 guesses\end{tabular}} & \multicolumn{1}{l|}{\begin{tabular}[c]{@{}l@{}}Max iteration\\ of OSQP=200\end{tabular}} & \multicolumn{1}{l|}{\begin{tabular}[c]{@{}l@{}}Max iteration\\ of OSQP=600\end{tabular}} & \multicolumn{1}{l}{\begin{tabular}[c]{@{}l@{}}Currently KNN\\ learner\end{tabular}} \\ \midrule
Success Rate          & 96.5\%               & 94.75\%              & 100 \%                                                                                 & 98.75\%                                                                                 & 1.25 \%                                                          & 0 \%                                                            & 78.25\%                                                         & 78.25\%                                                        & 99.25\%                                                                    & 98.5\%                                                                     & 92.5\%                                                                                                    & 94.5\%                                                                                  & 98.75\%                                                                                 & 96\%                                                                                \\ \midrule
Avg. Solve Time       & 260 ms               & 522 ms               & 790 ms                                                                                 & 616 ms                                                                                & 72 ms                                                            &                                                                 & 29 ms                                                           & 81 ms                                                          & 30 ms                                                                      & 96 ms                                                                      & 61 ms                                                                                                     & 65 ms                                                                                   & 80 ms                                                                                   & 16 ms                                                                               \\ \midrule
Max. Solve Time       & 1.29 sec             & 3.59 sec             & 17.3 s                                                                                 & 70 s                                                                                & 82 ms                                                            &                                                                 & 198 ms                                                          & 1.96 sec                                                       & 280 ms                                                                     & 990 ms                                                                     & 327 ms                                                                                                    & 190 ms                                                                                  & 400 ms                                                                                  & 174 ms                                                                              \\ \midrule
Avg. Objective        & -7250                & -7339                & -8606                                                                                  & -8435                                                                                 & -5703                                                            &                                                                 & -8620                                                           & -8952                                                          & -8258                                                                      & -8618                                                                      & -6461                                                                                                     & -8631                                                                                   & -8687                                                                                   & -8637                                                                               \\ \midrule
Avg. \# of iterations & 4.73                 & 4.6                  & 1                                                                                      & 1                                                                                     &                                                                  &                                                                 & 1                                                               & 1                                                              & 1.05                                                                       & 1.06                                                                       & 1.25                                                                                                      & 4.85                                                                                    & 3.85                                                                                    & 2.57                                                                                \\ \midrule
Solver                & Gurobi+Knitro        & Gurobi+IPOPT         & Gurobi                                                                                 & Gurobi                                                                                & Knitro                                                           & IPOPT                                                           & Knitro                                                          & IPOPT                                                          & Knitro                                                                     & IPOPT                                                                      & Knitro                                                                                                    & OSQP                                                                                    & OSQP                                                                                    & OSQP                                                                                \\ \bottomrule

\end{tabular}}
\label{Tab:benchmark}
\end{table*}

\subsection{Complementary Formulation}
An equivalent formulation of the bookshelf organization problem formulation is to turn all the integer variables $z_{i}$ into continuous variables with complementary constraints, such that the complete formulation can be solved through NLP solvers. We implemented the formulation from \cite{park2018semidefinite}:

\begin{equation}
    z_{i}(1-z_{i}) = 0
\end{equation}

Which is equivalent to the constraint $z_{i} \in \{0, 1\}$. We have also tried other implementations in \cite{stein2004continuous} and get similar or worse performances.

After collecting the dataset, we can learn the optimal problem-solution mapping from the problem parameter $\Theta$ to the solution $(\textbf{x}, \textbf{z})$. Online, the learner samples multiple warm-start points $(\textbf{x}, \textbf{z})$ for the NLP. The NLP with complementary constraints are then solved using the sampled warm-start one-by-one until a feasible solution is returned.

\subsection{MIP Formulation}
We convert bilinear constraints into mixed-integer linear constraints by gridding the solution space $\mathbb{S}$ and approximating the constraints locally inside grids with McCormick envelopes similar to \cite{dai2019global}. A McCormick envelope relaxation of one bilinear constraint $w=xy$ \cite{castro2015tightening} is the best linear approximation defined over a pair of lower and upper bounds $[x^{L}, x^{U}]$ and $[y^{L}, y^{U}]$. Therefore, we first assign grids $G(\textbf{x})$ to $\mathbb{S}$. Let $\{g_{l}\}$, $l=1,...,L$ be one cell in the grid with upper and lower bounds $[\textbf{x}^{L}, \textbf{x}^{U}]$. We introduce additional integer variables $\textbf{n}$, $n_{i} \in \{0, 1\}$. Each unique value of $\textbf{n}$ corresponds to one cell in the grid within which McCormick envelope relaxations are applied. The constraints $b_{j}(\textbf{x}, \textbf{z}; \Theta) \leq 0$, $j=1,...,m_{b}$ are converted into $L$ constraints: $E_{b_{j},l}(\textbf{x}, \textbf{z}, \textbf{n}; \Theta) \leq 0$ for $j = 1,...,m_{b}$ and $l = 1,...,L$, turning (\ref{Eqn:General_formulation}) into an MICP:

\begin{equation}
\begin{split}
& \underset{\textbf{x}, \ \textbf{z}, \ \textbf{n}}{\text{minimize}} \ f_{obj}(\textbf{x}, \textbf{z}, \textbf{n}; \Theta)\\
\text{s. t.} \qquad\; 
f_{i}(\textbf{x}, \textbf{z}; \Theta) \leq 0, & \ \ i = 1,...,m_{f} \\
E_{b_{j},l}(\textbf{x}, \textbf{z}, \textbf{n}; \Theta)& \leq 0, \ \ j = 1,...,m_{b}, \ \ l = 1,...,L \qquad
\label{Eqn:Updated_formulation}
\end{split}
\end{equation}

We use a $log_{2}N$ formulation which means $N$ grids are represented by $\ceil{log_{2}N}$ binary variables. For example, $17 \sim 32$ grids are represented by 5 integer variables. Since the intervals generated by the clustering methods are not necessarily connected, the proposed method in \cite{dai2019global} does not work. The formulation we use can incorporate intervals that are not connected hence more general. The price to pay is additional continuous variables. Please see Appendix for details of the formulation.  Despite the relatively smaller number of binary variables due to $log_{2}N$ formulation, approximating nonlinear constraints with mixed-integer convex constraints still generates a large number of integer variables when high approximation accuracy is desired. As a result, the formulation suffers from extended solving time. We use data-driven methods to reduce the number of integer variables and learn good warm-starts for online applications.

\subsubsection{Reduced Scale MIP Formulation}
\cite{lin2022reduce} proposes ReDUCE as a method to generate smaller scale MIP problems. The basic idea of using unsupervised learning to reduce the size of mixed-integer programs is to identify and cluster the important regions on the solution manifold, and generate MIP formulations with smaller scale. If $\mathbb{R}^{dim(\textbf{x})}$ is segmented into smaller regions, e.g. clusters, the required integer variables for each cluster can be reduced. This may be seen in Fig. \ref{Fig:Figure_1} which is an instance of 2-dimensional (dim) $X(\Theta)$ from a bookshelf experiment described in Sec.~\ref{Sec:experiment_setup}. 

This paragraph reviews the general steps of ReDUCE. To begin, ReDUCE uses pre-solved dataset $(\Theta_{k}, \textbf{x}_{k})$, $k=1,...,K$. We pre-assign grids $G(\textbf{x})$ to the solution space $\mathbb{S}$. The size of grids depends on the approximation accuracy requirement for bilinear constraints. ReDUCE begins by performing unsupervised learning on the pre-solved dataset to retrieve clusters that indicates regions on $X(\theta)$. Density-Based Spatial Clustering of Applications with Noise (DBSCAN) \cite{ester1996density} was used to cluster on $\textbf{x}$ which gives clusters ${\textbf{x}_{1}}, {\textbf{x}_{2}}, ..., {\textbf{x}_{C}}$, where $C$ is the number of clusters. We then trace $\Theta \rightarrow \textbf{x}$ map backwards to create clusters in $\Theta$ space, i.e. $(\{\Theta\}_{1}, \{\textbf{x}\}_{1}), ..., (\{\Theta\}_{C}, \{\textbf{x}\}_{C})$. Next, a supervised classifier is trained to classify $(\Theta, c)$. We use Random Forest which requires relatively smaller amounts of data to train than deep learning methods. For all $\{\textbf{x}\}_{c}$ in cluster $c$, we find the grid cells that they occupy and re-assign integer variables $\textbf{z}_{c}$ and $\textbf{n}_{c}$ to each cluster. At this point, the original MIP problem is segmented into smaller scale MIP problems with faster solving speed. When a new problem instance comes online, the classifier is used to point the new problem to one of the existing clusters and the reduced scale MIP is solved.

The idea of solving smaller sub-MIPs bears similarities with algorithms such as RENS \cite{berthold2014rens} and Neural Diving \cite{nair2020solving} where a subset of integer variables are fixed from linear program solutions or a learned model, while the others are solved. Our method, however, segments the problem through an unsupervised clustering approach.

\subsubsection{Convex Formulation}
In the extreme case, if the learner can provide a complete warm-start of the binary variables, the problem online becomes convex permitting fast QP solvers. \cite{bertsimas2021voice,cauligi2021coco,bertsimas2019online} define integer strategies to be tuples of $\mathcal{I}(\Theta_{i}) = (\textbf{z}^{*}, \textbf{n}^{*}, \mathcal{T}(\Theta_{i}))$, where $(\textbf{x}^{*}, \textbf{z}^{*}, \textbf{n}^{*})$ is an optimizer for problem \eqref{Eqn:Updated_formulation}, $\mathcal{T}(\Theta_{i}) = \{i \in {1,...,m_{f}}, \ l \in 1,...,L| f_{i}(\textbf{x}^{*}, \textbf{z}^{*}; \Theta) = 0, \ E_{b_{j},l}(\textbf{x}^{*}, \textbf{z}^{*}, \textbf{n}^{*}; \Theta)=0 \}$ is the set of active inequality constraints. Given an optimal integer strategy $\mathcal{I}(\Theta_{i})$, solutions to (\ref{Eqn:Updated_formulation}) can be retrieved through solving a convex optimization problem. This idea is adopted for learning warm-start with a convex formulation online. We have implemented both K-nearest neighbor to sample the points that are closest to the new problem $\Theta$. We have also tested a simple neural-network structure identical to \cite{cauligi2021coco}.

\subsection{ADMM Formulation}
We implemented a non-data-driven benchmark based on the ADMM algorithm in \cite{aydinoglu2022real}. The idea is to separate the MINLP into a MIP formulation and an NLP formulation with the exact copies of variables but different constraints, and iterate between them until the consensus is reached. The MIP formulation contains constraints $A, B, D, G, H1, I1, J1, K2, K3, L2, L3$, while missing bilinear constraints $C, E, F, K1, L1$. The NLP formulation contains constraints $A, B, C, D, E, F, H, I, J, K, L$, while missing the mode constraint G which says $z_{i}$s are integers. The formulation is elaborated in Eqn. \ref{eqn:MIP_NLP}.

\begin{equation}
\begin{aligned}
& \underset{\textbf{x}_{i}, \ \textbf{R}_{i}, \ \textbf{v}_{i,k}, \ \textbf{a}_{j}, \ b_{j}, \ z_{i}}{\text{minimize}} \quad f_{obj} \\
& \text{subject to} \\
& \text{\textbf{Mixed integer constraints:}} \\
& \text{A, B, D, G, H1, I1, J1, K2, K3, L2, L3} \\
& \text{\textbf{Nonlinear constraints:}} \\  
& \text{A, B, C, D, E, F, H, I, J, K, L}
\end{aligned}
\label{eqn:MIP_NLP}
\end{equation}

The ADMM iterates the following 3 steps: first it solves the MIP formulation in (\ref{eqn:MIP_NLP}), then solves the NLP formulation in (\ref{eqn:MIP_NLP}). Finally it updates the dual variables. The process continues till convergence. Let $\textbf{var}_{1}$ and $\textbf{var}_{2}$ be the copies of the variables in formulation (\ref{eqn:MIP_NLP}). Let $\textbf{w}$ be the dual variables. $\textbf{G}$ are the weights. The first step is:

\begin{equation*}
\begin{aligned}
& \underset{\textbf{var}_{1}}{\text{minimize}} \quad ||\textbf{var}_{1}^{i}-\textbf{var}_{2}^{i}+\textbf{w}^{i}||_{\textbf{G}_{k}} \\
& \text{s.t.  Mixed-integer constraints in ~\eqref{eqn:MIP_NLP}} 
\end{aligned}
\end{equation*}

The next step solves the NLP:
\begin{equation*}
\begin{aligned}
& \underset{\textbf{var}_{2}}{\text{minimize}} \quad ||\textbf{var}_{2}^{i}-(\textbf{var}_{1}^{i+1}+\textbf{w}^{i}))||_{\textbf{G}_{k}} \\
& \text{s.t.  Nonlinear constraints in ~\eqref{eqn:MIP_NLP}} 
\end{aligned}
\end{equation*}

Finally, we update the dual variables and weights:

\begin{equation}
\begin{aligned}
    \textbf{w}^{i+1} = \textbf{w}^{i} + \textbf{var}_{1}^{i} - \textbf{var}_{2}^{i} \\
    \textbf{G}_{k+1} = \gamma\textbf{G}_{k} \\
    \textbf{w}^{i+1} = \textbf{w}^{i} / \gamma
\end{aligned}
\end{equation}

\section{Experiment} 
\subsection{Setup}
\label{Sec:experiment_setup}
We place 3 books inside the shelf where 1 additional book is to be inserted. Grids are assigned to the variables involved in the non-convex constraint \textbf{C}, \textbf{E}, \textbf{F}, \textbf{K1}, \textbf{L1}: $\textbf{R}_{i}(\theta_{i})$, $\textbf{a}_{j}$ and $\textbf{v}_{i,k}$. These variables span a 48-dim space. The rotation angles $\theta_{i}$, which includes $\textbf{R}_{i}$, are gridded at $\frac{\pi}{8}$ intervals. Elements in $\textbf{a}_{j}$ are gridded on 0.25 intervals. Elements in $\textbf{v}_{i,k}$ are gridded at intervals $\frac{1}{4}$ the shelf width $W$ and height $H$. Table \ref{Tab:variable_range} lists the number of grids for each non-convex variable. Our MICP formulation results in 130 integer variables in total. The input vector includes the center positions, angles, heights and widths of stored books and height and width of the book to be inserted. The input dimension is 17.

% How the data was collected/generated
The problem instances are generated using a 2-dim simulated environment of books on a shelf. Initially, 4 randomly sized books are arbitrarily placed on the shelf, and then 1 is randomly removed and regarded as the book to be inserted. Contrary to the sequence, the initial state with 4 books represents one feasible (not necessarily optimal) solution to the problem of placing a book on a shelf with 3 existing books. This guarantees that all problem instances are feasible. Since this problem can be viewed as high-level planning for robotic systems, the simulated data is sufficient. For applications outside of the scope of this paper real-world data may be preferable in this pipeline.

All methods are tested on 400 randomly sampled bookshelf problems in the same distribution using the same method mentioned above. All results are listed in table \ref{Tab:benchmark}.

\subsection{Complementary Formulation}
\label{Sec:Complementary_formulation}
As mentioned in the Section \ref{Sec:Complementary_formulation}, we turn the binary variables into continuous variables with complementary constraints. We selected both IPOPT and Knitro solvers to solve the MPCC formulation. A small relaxation of $\epsilon = 10^{-3}$ on the right hand side of the complementary constraints is used to increase the chance of getting a feasible solution.

\textit{Zero Warm Start}
As a baseline, we solve the problem with the default initial guess of the solvers (all zeros). The results are shown in the column of Table \ref{Tab:benchmark} named "MPCC\_default". As we can see, despite the relaxation of complementary constraints, the rate of solving the problem successfully with zero warm-start is almost zero.

\textit{Warm Start with Manually Designed Heuristics}
Another baseline is to warm-start the formulation with human-designed heuristics. The bookshelf problem requires inserting one book minimizing the movements of the existing books. Therefore, a good initial guess is to directly use the continuous and discrete variables from the original scene. We also pre-solve the separating planes using the original scene to warm-start the parameters $\textbf{a}$ and $b$. The authors believe this is the best initial guess readily available without paying much more effort. The results are shown in Table \ref{Tab:benchmark} in the column "MPCC\_Manual".

Implementing this initial guess increases the success rate to more than $70\%$, significantly better than the default initial guess from the solver. However, we noticed that this success rate drops as the number of books increases. When there are only 2 books on the shelf, inserting the 3rd book has over $80\%$ success rate with the same approach for the initial guess. Therefore, the manually selected initial guess may not scale to problems with much larger sizes.

\textit{Warm Start with pre-solved Dataset}
We show the result of using K-nearest neighbor learner to pick out the top 3 candidates from a 15000 pre-solved dataset in the column "MPCC\_KNN3" column in Table \ref{Tab:benchmark}. The results show that a simple KNN solver is able to significantly increase the rate of getting a feasible solution to more than $98\%$. The average iteration of solving the problem is only slightly more than 1 showing that most problem instances are solved with only 1 K-nearest neighbor sample. Since the MPCC formulation keeps the non-convexity of the original formulation, it permits the solver to explore larger regions, hence the learner does not need to pick out a point or region that is very close to the optimal solution as required by the convex formulation. Note that we found warm-starting the non-convex continuous variables is as important as warm-starting the discrete (complementary) variables.

% If z only, the success rate drop to 67\%, -8168, avg time 0.256s, max time 0.798s, iteration 1.466.

\subsection{Mixed-integer Programming Formulation}
\label{Sec:MIP_formulation}
\subsubsection{Reduced Scale MIP Formulation}
We implemented the clustering method proposed in section \ref{Sec:MIP_formulation} to reduce the size of MIP formulation. The dataset used for clustering has a size of 1000, and is classified into 26 clusters. A random forest classifier is trained to predict the cluster label on the dataset which achieves 92 \%. As the new problem instance arrives, the random forest classifier selects the cluster, then the region occupied by the non-convex solutions for this dataset is used to solve the mixed-integer programming. The results are shown in the columns named "MIP" in Table \ref{Tab:benchmark}.

As the results show, this method has the highest chance to solve the new problem instances. However, the solving time is significantly longer than the MPCC formulation and the convex formulation. The reason is that some problems which do not get classified correctly. It can take minutes to solve those cases making the average solving time much longer. We are still working on tuning the clustering method to eliminate those cases.

\subsubsection{Convex Formulation}
To make the online formulation convex, the non-convex part of the continuous solutions is first discretized, then the complete list of integer variables is used for training. We discretize the variables using the range given by Table \ref{Tab:variable_range}. A K-nearest neighbor learner is trained to pick out the first 10 integer solutions. They are used to turn the mixed-integer program into convex programs which are then solved one by one until a feasible solution is retrieved. We select the solver OSQP to solve those convex problems. This solver automatically uses the previous solution to warm-start the next solution for faster solving speed. The maximal iterations can be tuned. As most of the problems take no more than 200 iterations to solve, we set this parameter to 200. Setting the maximal iteration to 600 will incorporate more solutions, but the average solving time is longer as it takes more time before the solver decides infeasibility. The result is shown in the columns named "Convex" in Table \ref{Tab:benchmark}.

The result shows that the convex formulation achieves the fastest solving speed due to the strong performance of convex solvers. The authors suspect that the solvers such as OSQP utilized the heuristics specific to convexity, hence run faster than well-initialized NLP solvers such as Knitro. On the other hand, the price to pay for turning the problem into convex is that the learner needs to perform better than in the MPCC case. This is indicated by the 4.85 average iterations required to find a feasible solution. Due to the inherent issue with KNN, some instances are classified incorrectly, and the solver wastes much longer time going over the infeasible integer variables. 

\begin{table}[]
\centering
\caption{Segmentations of Non-convex Variables}
\begin{tabular}{l|c|c}
\hline
\multicolumn{1}{c|}{variable}   & range             & \# of segmentations \\ \hline
Item orientation $\theta$ (rad) & {[}-pi/2, pi/2{]} & 8                   \\ \hline
Separating plane normal $a$     & {[}-1, 1{]}       & 8                   \\ \hline
Vertex x position $v_{x}$ (cm)  & {[}-9, 9{]}       & 4                   \\ \hline
Vertex y position $v_{y}$ (cm)  & {[}0, 11{]}       & 4                   \\ \hline
\end{tabular}
\label{Tab:variable_range}
\end{table}

As the space is discretized, the $\Theta$ space is no longer Euclidean. The quantities that are close to each other in the Euclidean distance sense, but on different sides of the discretization boundary will be assigned different integer variables, hence no longer close to each other. The KNN utilizes Euclidean distance, thus ignoring this non-Euclidean effect and failing to select the correct warm-start integers in those cases.

\subsection{Nonlinear ADMM Formulation}
\label{Sec:ADMM_formulation}
As a non-data-driven benchmark, we implemented the nonlinear ADMM method on the MINLP formulation by iterating between a mixed-integer formulation and a nonlinear formulation as described in Section \ref{Sec:ADMM_formulation}. To keep this method non-data-driven, manually designed heuristics same as \ref{Sec:Complementary_formulation} are used to warm-start the nonlinear formulation to increase the feasibility rate. A sufficient effort is devoted to tuning the weights $\textbf{G}$, $\gamma$, and scaling the variables properly to ensure the performance is on the better end. Since the mixed-integer approximations on the bilinear constraints induce numerical errors, termination conditions and accuracy of bilinear collision avoidance constraints are set lower to match the accuracy of mixed-integer envelope formulation. The results are shown in the column named "ADMM" in Table \ref{Tab:benchmark}.

According to the results, ADMM works fairly well for getting a feasible solution. The solving speed is several times slower compared to the data-driven MPCC formulation and convex formulation. One important reason is that it usually takes 4-5 iterations (1 iteration includes 1 MIP and 1 NLP) to get a feasible solution. The nonlinear formulation takes more than 70\% of the solving time, due to the separating plane constraints. If data-driven methods can be used to warm-start the nonlinear formulation, a faster solving speed may be achieved. In addition, the average optimal cost from ADMM is worse than most of the data-driven methods. The authors suspect that the reason is due to the extra consensus terms in the objective function which guides the convergence.

\subsection{Hardware Experiment}
\label{Sec:hardware_experiment}
We implement our optimization results on hardware with a 6 degree of freedom manipulator \cite{noh2020minimal} to insert a book onto a shelf. To automate the system, a trajectory planner is required to generate the insertion motion, which is common in practice. As the focus of this paper is on high level planning, we simplify this part by manually selecting waypoints. The hardware implementation is shown in the attached video.

% We implement several optimization results on hardware by commanding a 6 degree of freedom manipulator to insert 1 book onto a shelf. To automate the complete system, a trajectory planner is required to generate motion for insertion, which is achievable given vast literature in this field. As the focus of this paper is to plan the final scene, we simplify this part by manually selecting waypoints. The success of the hardware implementation is shown in the video attached to this paper.

\section{Conclusion, Discussion and Future Work} 
\label{Sec:conclusion}
This paper proposes data-driven methods to solve MPCC, MIP, or convex formulation applied to a bookshelf organization problem formulated into mixed-integer non-convex optimization incorporating mode switch and separating planes as collision avoidance constraints. The performances of the proposed data-driven methods are compared. The planned organizations are implemented on the hardware.

According to our test, data-driven methods can simultaneously achieve good optimality, fast solving speed, and high success rate leveraging on a reasonable amount of pre-solved data. Among the online formulations, MPCC keeps the original non-convex formulations hence can converge better even from a relatively bad initial warm-start. The implementation of MPCC on NLP solvers is easier than the convex formulation which requires converting the non-convex constraints into mixed-integer envelope constraints. The solving speed of MPCC formulation is also fast. The convex formulation has the fastest solving speed due to well-performed convex solvers. However, since additional constraints are added to make the problem convex, it is more difficult for the learner to select a good warm-start, resulting in more trials. Finally, non-data-driven methods such as nonlinear ADMM have a decent chance to get a feasible solution. However, both the optimality and the solving speed are much worse.

The obvious merit of the data-driven approach is that it can solve the problem fast, but also optimally and reliably as long as the new problem instance is within the trained parameter regime. An obvious challenge with this method is how to generalize the training to cases outside the dataset. To achieve this goal, simple learners such as KNN do not perform well, and more modern learning methods such as auto-encoders are to be explored.

% \section*{Acknowledgements}
\textbf{\textit{Acknowledgements}} Special thanks to Abhishek Cauligi and Professor Bartolomeo Stellato for helpful discussions, and Yusuke Tanaka and Hyunwoo Nam for hardware help.

% % Appendixes should appear before the acknowledgment.

{
\bibliographystyle{IEEEtran}
\bibliography{references}

\begin{thebibliography}{10}
\providecommand{\url}[1]{#1}
\csname url@samestyle\endcsname
\providecommand{\newblock}{\relax}
\providecommand{\bibinfo}[2]{#2}
\providecommand{\BIBentrySTDinterwordspacing}{\spaceskip=0pt\relax}
\providecommand{\BIBentryALTinterwordstretchfactor}{4}
\providecommand{\BIBentryALTinterwordspacing}{\spaceskip=\fontdimen2\font plus
\BIBentryALTinterwordstretchfactor\fontdimen3\font minus
  \fontdimen4\font\relax}
\providecommand{\BIBforeignlanguage}[2]{{%
\expandafter\ifx\csname l@#1\endcsname\relax
\typeout{** WARNING: IEEEtran.bst: No hyphenation pattern has been}%
\typeout{** loaded for the language `#1'. Using the pattern for}%
\typeout{** the default language instead.}%
\else
\language=\csname l@#1\endcsname
\fi
#2}}
\providecommand{\BIBdecl}{\relax}
\BIBdecl

\bibitem{deits2014footstep}
R.~Deits and R.~Tedrake, ``Footstep planning on uneven terrain with
  mixed-integer convex optimization,'' in \emph{2014 IEEE-RAS international
  conference on humanoid robots}.\hskip 1em plus 0.5em minus 0.4em\relax IEEE,
  2014, pp. 279--286.

\bibitem{lin2019optimization}
X.~Lin, J.~Zhang, J.~Shen, G.~Fernandez, and D.~W. Hong, ``Optimization based
  motion planning for multi-limbed vertical climbing robots,'' in \emph{2019
  IEEE/RSJ International Conference on Intelligent Robots and Systems
  (IROS)}.\hskip 1em plus 0.5em minus 0.4em\relax IEEE, 2019, pp. 1918--1925.

\bibitem{dai2014whole}
H.~Dai, A.~Valenzuela, and R.~Tedrake, ``Whole-body motion planning with
  centroidal dynamics and full kinematics,'' in \emph{2014 IEEE-RAS
  International Conference on Humanoid Robots}.\hskip 1em plus 0.5em minus
  0.4em\relax IEEE, 2014, pp. 295--302.

\bibitem{winkler2018gait}
A.~W. Winkler, C.~D. Bellicoso, M.~Hutter, and J.~Buchli, ``Gait and trajectory
  optimization for legged systems through phase-based end-effector
  parameterization,'' \emph{IEEE Robotics and Automation Letters}, vol.~3,
  no.~3, pp. 1560--1567, 2018.

\bibitem{shirai2020risk}
Y.~Shirai, X.~Lin, Y.~Tanaka, A.~Mehta, and D.~Hong, ``Risk-aware motion
  planning for a limbed robot with stochastic gripping forces using nonlinear
  programming,'' \emph{IEEE Robotics and Automation Letters}, vol.~5, no.~4,
  pp. 4994--5001, 2020.

\bibitem{soler2011route}
M.~Soler, A.~Olivares, P.~Bonami, and E.~Staffetti, ``En-route optimal flight
  planning constrained to pass through waypoints using minlp,'' 2011.

\bibitem{boyd2007branch}
S.~Boyd and J.~Mattingley, ``Branch and bound methods,'' \emph{Notes for
  EE364b, Stanford University}, pp. 2006--07, 2007.

\bibitem{tordesillas2019faster}
J.~Tordesillas, B.~T. Lopez, and J.~P. How, ``Faster: Fast and safe trajectory
  planner for flights in unknown environments,'' in \emph{2019 IEEE/RSJ
  international conference on intelligent robots and systems (IROS)}.\hskip 1em
  plus 0.5em minus 0.4em\relax IEEE, 2019, pp. 1934--1940.

\bibitem{lin2021designing}
X.~Lin, M.~S. Ahn, and D.~Hong, ``Designing multi-stage coupled convex
  programming with data-driven mccormick envelope relaxations for motion
  planning,'' in \emph{2021 IEEE/RSJ International Conference on Robotics and
  Automation (ICRA)}.\hskip 1em plus 0.5em minus 0.4em\relax IEEE, 2021.

\bibitem{nair2020solving}
V.~Nair, S.~Bartunov, F.~Gimeno, I.~von Glehn, P.~Lichocki, I.~Lobov,
  B.~O'Donoghue, N.~Sonnerat, C.~Tjandraatmadja, P.~Wang \emph{et~al.},
  ``Solving mixed integer programs using neural networks,'' \emph{arXiv
  preprint arXiv:2012.13349}, 2020.

\bibitem{tang2020reinforcement}
Y.~Tang, S.~Agrawal, and Y.~Faenza, ``Reinforcement learning for integer
  programming: Learning to cut,'' in \emph{International Conference on Machine
  Learning}.\hskip 1em plus 0.5em minus 0.4em\relax PMLR, 2020, pp. 9367--9376.

\bibitem{zhu2019fast}
J.-J. Zhu and G.~Martius, ``Fast non-parametric learning to accelerate
  mixed-integer programming for online hybrid model predictive control,''
  \emph{arXiv preprint arXiv:1911.09214}, 2019.

\bibitem{cauligi2021coco}
A.~Cauligi, P.~Culbertson, E.~Schmerling, M.~Schwager, B.~Stellato, and
  M.~Pavone, ``Coco: Online mixed-integer control via supervised learning,''
  \emph{arXiv preprint arXiv:2107.08143}, 2021.

\bibitem{deits2015computing}
R.~Deits and R.~Tedrake, ``Computing large convex regions of obstacle-free
  space through semidefinite programming,'' in \emph{Algorithmic foundations of
  robotics XI}.\hskip 1em plus 0.5em minus 0.4em\relax Springer, 2015, pp.
  109--124.

\bibitem{dai2017distributed}
L.~Dai, Q.~Cao, Y.~Xia, and Y.~Gao, ``Distributed mpc for formation of
  multi-agent systems with collision avoidance and obstacle avoidance,''
  \emph{Journal of the Franklin Institute}, vol. 354, no.~4, pp. 2068--2085,
  2017.

\bibitem{raghunathan2022pyrobocop}
A.~U. Raghunathan, D.~K. Jha, and D.~Romeres, ``Pyrobocop: Python-based robotic
  control \& optimization package for manipulation,'' in \emph{2022
  International Conference on Robotics and Automation (ICRA)}.\hskip 1em plus
  0.5em minus 0.4em\relax IEEE, 2022, pp. 985--991.

\bibitem{stolle2006policies}
M.~Stolle and C.~G. Atkeson, ``Policies based on trajectory libraries,'' in
  \emph{Proceedings 2006 IEEE International Conference on Robotics and
  Automation, 2006. ICRA 2006.}\hskip 1em plus 0.5em minus 0.4em\relax IEEE,
  2006, pp. 3344--3349.

\bibitem{liu2013biped}
C.~Liu, C.~G. Atkeson, and J.~Su, ``Biped walking control using a trajectory
  library,'' \emph{Robotica}, vol.~31, no.~2, pp. 311--322, 2013.

\bibitem{tassa2007receding}
Y.~Tassa, T.~Erez, and W.~Smart, ``Receding horizon differential dynamic
  programming,'' \emph{Advances in neural information processing systems},
  vol.~20, 2007.

\bibitem{fiacco1983introduction}
A.~V. Fiacco, \emph{Introduction to sensitivity and stability analysis in
  nonlinear programming}.\hskip 1em plus 0.5em minus 0.4em\relax Academic
  press, 1983, vol. 165.

\bibitem{pistikopoulos2002line}
E.~N. Pistikopoulos, V.~Dua, N.~A. Bozinis, A.~Bemporad, and M.~Morari,
  ``On-line optimization via off-line parametric optimization tools,''
  \emph{Computers \& Chemical Engineering}, vol.~26, no.~2, pp. 175--185, 2002.

\bibitem{gal2010postoptimal}
T.~Gal, \emph{Postoptimal Analyses, Parametric Programming, and Related Topics:
  degeneracy, multicriteria decision making, redundancy}.\hskip 1em plus 0.5em
  minus 0.4em\relax Walter de Gruyter, 2010.

\bibitem{bemporad2002explicit}
A.~Bemporad, M.~Morari, V.~Dua, and E.~N. Pistikopoulos, ``The explicit linear
  quadratic regulator for constrained systems,'' \emph{Automatica}, vol.~38,
  no.~1, pp. 3--20, 2002.

\bibitem{dua1999algorithms}
V.~Dua and E.~N. Pistikopoulos, ``Algorithms for the solution of
  multiparametric mixed-integer nonlinear optimization problems,''
  \emph{Industrial \& engineering chemistry research}, vol.~38, no.~10, pp.
  3976--3987, 1999.

\bibitem{tondel2003algorithm}
P.~T{\o}ndel, T.~A. Johansen, and A.~Bemporad, ``An algorithm for
  multi-parametric quadratic programming and explicit mpc solutions,''
  \emph{Automatica}, vol.~39, no.~3, pp. 489--497, 2003.

\bibitem{dua2000algorithm}
V.~Dua and E.~N. Pistikopoulos, ``An algorithm for the solution of
  multiparametric mixed integer linear programming problems,'' \emph{Annals of
  operations research}, vol.~99, no.~1, pp. 123--139, 2000.

\bibitem{hauser2016learning}
K.~Hauser, ``Learning the problem-optimum map: Analysis and application to
  global optimization in robotics,'' \emph{IEEE Transactions on Robotics},
  vol.~33, no.~1, pp. 141--152, 2016.

\bibitem{tang2019data}
G.~Tang and K.~Hauser, ``A data-driven indirect method for nonlinear optimal
  control,'' \emph{Astrodynamics}, vol.~3, no.~4, pp. 345--359, 2019.

\bibitem{power2022variational}
T.~Power and D.~Berenson, ``Variational inference mpc using normalizing flows
  and out-of-distribution projection,'' \emph{arXiv preprint arXiv:2205.04667},
  2022.

\bibitem{dai2019global}
H.~Dai, G.~Izatt, and R.~Tedrake, ``Global inverse kinematics via mixed-integer
  convex optimization,'' \emph{The International Journal of Robotics Research},
  vol.~38, no. 12-13, pp. 1420--1441, 2019.

\bibitem{marchand2003cutting}
H.~Marchand, A.~Martin, R.~Weismantel, and L.~Wolsey, ``Cutting planes in
  integer and mixed integer,'' \emph{Discrete Optimization: The State of the
  Art}, vol.~11, p. 397, 2003.

\bibitem{williams2013model}
H.~P. Williams, \emph{Model building in mathematical programming}.\hskip 1em
  plus 0.5em minus 0.4em\relax John Wiley \& Sons, 2013.

\bibitem{anitescu1997formulating}
M.~Anitescu and F.~A. Potra, ``Formulating dynamic multi-rigid-body contact
  problems with friction as solvable linear complementarity problems,''
  \emph{Nonlinear Dynamics}, vol.~14, no.~3, pp. 231--247, 1997.

\bibitem{drumwright2015rapidly}
E.~Drumwright, ``Rapidly computable viscous friction and no-slip rigid contact
  models,'' \emph{arXiv preprint arXiv:1504.00719}, 2015.

\bibitem{kojima1991unified}
M.~Kojima, N.~Megiddo, T.~Noma, and A.~Yoshise, ``A unified approach to
  interior point algorithms for linear complementarity problems: A summary,''
  \emph{Operations Research Letters}, vol.~10, no.~5, pp. 247--254, 1991.

\bibitem{posa2014direct}
M.~Posa, C.~Cantu, and R.~Tedrake, ``A direct method for trajectory
  optimization of rigid bodies through contact,'' \emph{The International
  Journal of Robotics Research}, vol.~33, no.~1, pp. 69--81, 2014.

\bibitem{zhang2021transition}
J.~Zhang, X.~Lin, and D.~W. Hong, ``Transition motion planning for multi-limbed
  vertical climbing robots using complementarity constraints,'' in \emph{2021
  IEEE International Conference on Robotics and Automation (ICRA)}.\hskip 1em
  plus 0.5em minus 0.4em\relax IEEE, 2021, pp. 2033--2039.

\bibitem{cleac2021fast}
S.~L. Cleac'h, T.~Howell, M.~Schwager, and Z.~Manchester, ``Fast
  contact-implicit model-predictive control,'' \emph{arXiv preprint
  arXiv:2107.05616}, 2021.

\bibitem{park2018semidefinite}
J.~Park and S.~Boyd, ``A semidefinite programming method for integer convex
  quadratic minimization,'' \emph{Optimization Letters}, vol.~12, no.~3, pp.
  499--518, 2018.

\bibitem{boyd2004convex}
S.~Boyd, S.~P. Boyd, and L.~Vandenberghe, \emph{Convex optimization}.\hskip 1em
  plus 0.5em minus 0.4em\relax Cambridge university press, 2004.

\bibitem{vielma2015mixed}
J.~P. Vielma, ``Mixed integer linear programming formulation techniques,''
  \emph{Siam Review}, vol.~57, no.~1, pp. 3--57, 2015.

\bibitem{posa2015stability}
M.~Posa, M.~Tobenkin, and R.~Tedrake, ``Stability analysis and control of
  rigid-body systems with impacts and friction,'' \emph{IEEE Transactions on
  Automatic Control}, vol.~61, no.~6, pp. 1423--1437, 2015.

\bibitem{bonami2007bonmin}
P.~Bonami and J.~Lee, ``Bonmin user’s manual,'' \emph{Numer Math}, vol.~4,
  pp. 1--32, 2007.

\bibitem{stellato2020osqp}
B.~Stellato, G.~Banjac, P.~Goulart, A.~Bemporad, and S.~Boyd, ``Osqp: An
  operator splitting solver for quadratic programs,'' \emph{Mathematical
  Programming Computation}, vol.~12, no.~4, pp. 637--672, 2020.

\bibitem{risk_aware}
Y.~Shirai \emph{et~al.}, ``Simultaneous contact-rich grasping and locomotion
  via distributed optimization enabling free-climbing for multi-limbed
  robots,'' \emph{IEEE Proc. 2022 IEEE/RSJ Int. Conf. Intell. Rob. Syst}, 2022.

\bibitem{stein2004continuous}
O.~Stein, J.~Oldenburg, and W.~Marquardt, ``Continuous reformulations of
  discrete--continuous optimization problems,'' \emph{Computers \& chemical
  engineering}, vol.~28, no.~10, pp. 1951--1966, 2004.

\bibitem{castro2015tightening}
P.~M. Castro, ``Tightening piecewise mccormick relaxations for bilinear
  problems,'' \emph{Computers \& Chemical Engineering}, vol.~72, pp. 300--311,
  2015.

\bibitem{lin2022reduce}
X.~Lin, G.~I. Fernandez, and D.~W. Hong, ``Reduce: Reformulation of mixed
  integer programs using data from unsupervised clusters for learning efficient
  strategies,'' in \emph{2022 International Conference on Robotics and
  Automation (ICRA)}.\hskip 1em plus 0.5em minus 0.4em\relax IEEE, 2022, pp.
  4459--4465.

\bibitem{ester1996density}
M.~Ester, H.-P. Kriegel, J.~Sander, and X.~Xu, ``Density-based spatial
  clustering of applications with noise,'' in \emph{Int. Conf. Knowledge
  Discovery and Data Mining}, vol. 240, 1996, p.~6.

\bibitem{berthold2014rens}
T.~Berthold, ``Rens,'' \emph{Mathematical Programming Computation}, vol.~6,
  no.~1, pp. 33--54, 2014.

\bibitem{bertsimas2021voice}
D.~Bertsimas and B.~Stellato, ``The voice of optimization,'' \emph{Machine
  Learning}, vol. 110, no.~2, pp. 249--277, 2021.

\bibitem{bertsimas2019online}
------, ``Online mixed-integer optimization in milliseconds,'' \emph{arXiv
  preprint arXiv:1907.02206}, 2019.

\bibitem{aydinoglu2022real}
A.~Aydinoglu and M.~Posa, ``Real-time multi-contact model predictive control
  via admm,'' in \emph{2022 International Conference on Robotics and Automation
  (ICRA)}.\hskip 1em plus 0.5em minus 0.4em\relax IEEE, 2022, pp. 3414--3421.

\bibitem{noh2020minimal}
D.~Noh, Y.~Liu, F.~Rafeedi, H.~Nam, K.~Gillespie, J.-s. Yi, T.~Zhu, Q.~Xu, and
  D.~Hong, ``Minimal degree of freedom dual-arm manipulation platform with
  coupling body joint for diverse cooking tasks,'' in \emph{2020 17th
  International Conference on Ubiquitous Robots (UR)}.\hskip 1em plus 0.5em
  minus 0.4em\relax IEEE, 2020, pp. 225--232.

\bibitem{belotti2011disjunctive}
P.~Belotti, L.~Liberti, A.~Lodi, G.~Nannicini, A.~Tramontani \emph{et~al.},
  ``Disjunctive inequalities: applications and extensions,'' \emph{Wiley
  Encyclopedia of Operations Research and Management Science}, vol.~2, pp.
  1441--1450, 2011.

\bibitem{vielma2011modeling}
J.~P. Vielma and G.~L. Nemhauser, ``Modeling disjunctive constraints with a
  logarithmic number of binary variables and constraints,'' \emph{Mathematical
  Programming}, vol. 128, no.~1, pp. 49--72, 2011.

\bibitem{avis1992pivoting}
D.~Avis and K.~Fukuda, ``A pivoting algorithm for convex hulls and vertex
  enumeration of arrangements and polyhedra,'' \emph{Discrete \& Computational
  Geometry}, vol.~8, no.~3, pp. 295--313, 1992.

\end{thebibliography}
}

\appendix
We explain the formulation that is used to solve the MIPs in Sec.~\ref{Sec:experiment_setup}.

For a non-convex constraint, we segment it into multiple regions and locally approximate or relax them into convex constraints. In this paper, we relax the bilinear constraints locally into convex polytopes (McCormick envelopes). Each polytope is associated with a unique combination of integer variable values, hence, mixed-integer convex constraints. Assume the number of regions used is $N$. Depending on the number of integer variables used, the formulation can generally be divided into 2 categories: 1) If the number of integer variables is $N$, we call it \textit{$N$ formulation}, e.g., the convex hull formulation \cite{belotti2011disjunctive}. 2) If the number of integer variables is $log_{2}N$, we call it \textit{$log_{2}N$ formulation}, e.g., \cite{vielma2011modeling}. \cite{vielma2011modeling} presents a \textit{$log_{2}N$ formulation} to model the special ordered sets of type 2 (sos2). However, there are several limitations of this formulation. For example, the segmented regions need to be connected to be a valid sos2 constraint, i.e., the two consecutive couple of non-zero entries can be any consecutive couple in the set. In this paper, we use MIP to model convex polytopes of arbitrary locations. This can increase the complexity for sos2 techniques as the polytopes can be disjunctive. The disjunctive constraints can be handled with convex hull formulations at a price of introducing more integer variables which may result in slower solving speeds.
% change difficulty to complexity--two sentences above

In this appendix, we demonstrate an intuitive but general \textit{$log_{2}N$ formulation} to model combinations of convex polytopes at any locations which serves as the base MIP formulation for the bilinear constraints in our paper. Assume the variable $\textbf{x}$ is enforced to be within one of the $N$ convex polytopes, denoted by $\textbf{A}_{i}\textbf{x} \leq \textbf{b}_{i}$, $i=1,...,N$. We introduce $m=log_{2}N$ binary variables $z_{1}, ..., z_{m}$, $z_{i} \in \{0,1\}$. Each combination of unique values of binary variables can be assigned to a convex polytope. Let the assignment be:
% Assume or we just assign them as such then simplify

\begin{equation}
    \textbf{z}=\Bar{\textbf{z}}_{i} \Rightarrow \textbf{A}_{i}\textbf{x} \leq \textbf{b}_{i}
\end{equation}

Where $\Bar{\textbf{z}}_{i} = [\Bar{z}_{i,1}, ..., \Bar{z}_{i,m}]$ are constant binary values associated with polytope $i$. Note $\Bar{z}_{i} \neq \Bar{z}_{j}$ if $i \neq j$. In other words, we require that when $\textbf{z}=\Bar{\textbf{z}}_{i}$, $\textbf{x}$ stays within the polytope $\textbf{A}_{i}\textbf{x} \leq \textbf{b}_{i}$; otherwise, the constraint is unenforced.

Denote the vertices of the polytopes by $\textbf{v}_{i,1}, ..., \textbf{v}_{i,n_{i}}$, $i=1, ..., N$, where $n_{i}$ is the number of vertices associated with polytope $i$. Each vertex $\textbf{v}_{i,j}$ is assigned a continuous non-negative variable $\lambda_{i,j} \in [0, 1]$. In general, one can run a mathematical program (e.g. \cite{avis1992pivoting}) to get vertices from the nondegenerate system of inequalities $\textbf{A}_{i}\textbf{x} \leq \textbf{b}_{i}$. As a result, the assignment becomes: 

\begin{equation}
    \textbf{z}=\Bar{\textbf{z}}_{i} \Rightarrow \ \ 
\begin{aligned}
    & \textbf{x} = \sum_{j=1}^{n_{i}} \lambda_{i,j} \textbf{v}_{i,j} \\
    & \sum_{j=1}^{n_{i}} \lambda_{i,j} = 1, \ \ \lambda_{i,j} \in [0, 1]
\end{aligned}
\label{Eqn:MICP_formulation_one}
\end{equation}

The formulation can be written as:

\begin{equation}
\begin{aligned}
    & (5.a) \quad \textbf{x} = \sum_{i=1}^{N} \sum_{j=1}^{n_{i}} \lambda_{i,j} \textbf{v}_{i,j} \\
    & (5.b) \quad  \sum_{i=1}^{N} \sum_{j=1}^{n_{i}} \lambda_{i,j} = 1, \ \ \lambda_{i,j} \in [0, 1] \\
    & (5.c) \quad  \sum_{k=1,...,N}^{k \neq i} \sum_{j=1}^{n_{k}} \lambda_{k,j} \leq \sum_{l=1}^{m} |z_{l} - \Bar{z}_{i,l}|
\end{aligned}
\label{Eqn:MICP_formulation_all}
\end{equation}

The set of constraint in \eqref{Eqn:MICP_formulation_all} enforces \eqref{Eqn:MICP_formulation_one} for all polytopes, as $\sum_{l=1}^{m} |z_{l} - \Bar{z}_{i,l}| = 0$ only when $\textbf{z}=\Bar{\textbf{z}}_{i}$, enforcing that all $\lambda$'s that are not associated with polytope $i$ to be zero. If $\textbf{z} \neq \Bar{\textbf{z}}_{i}$, $\sum_{l=1}^{m} |z_{l} - \Bar{z}_{i,l}| \geq 1$ and constraint (5.c) is looser than constraint (5.b), hence, trivial.

Note that formulation \eqref{Eqn:MICP_formulation_all} works for any convex polytope that can be written as $\textbf{A}_{i}\textbf{x}\leq \textbf{b}_{i}$. In this paper, the polytopes are McCormick envelope constraints which is a special case.

\end{document}